\documentclass[12pt,leqno]{article}

\begin{document}

\title{Potpourri, 5}

\author{Stephen William Semmes	\\
	Rice University		\\
	Houston, Texas}

\date{}

\maketitle

%%%%%%   a footnote for the first page 

\renewcommand{\thefootnote}{}   %%%%  this is to get rid of the unneeded
                                %%%%  marker for the following footnote
                                %%%%  (which is the only footnote in the
                                %%%%  document).

\footnotetext{These notes are connected to a ``potpourri'' class in
the mathematics department at Rice University.}

	Let $\{a_j\}_{j=1}^\infty$ be a sequence of nonnegative real
numbers which is \emph{submultiplicative} in the sense that
\begin{equation}
	a_{j + l} \le a_j \, a_l
\end{equation}
for all positive integers $j$, $l$.  Let us show that the limit
\begin{equation}
	\lim_{l \to \infty} (a_l)^{1/l}
\end{equation}
exists, and in fact is equal to the infimum of the $a_j$'s.  If $a_j =
0$ for some $j$, then $a_l = 0$ for all $l \ge j$, and so we may as
well assume that $a_j > 0$ for all $j$.

	Notice first that
\begin{equation}
	a_j \le (a_1)^j
\end{equation}
for all positive integers $j$, so that the sequence consisting of
$(a_l)^{1/l}$ for all positive integers $l$ is bounded.  More
generally, $a_{p \, k} \le (a_k)^p$ for all positive integers $k$ and
$p$, so that $(a_n)^{1/n} \le (a_k)^{1/k}$ when $n$ is a multiple of
$k$.

	Fix a positive integer $k$, and let $n$ be a positive integer
such that $k \le n$.  We can write $n$ as $p \, k + r$, where $p$, $r$
are integers, $p \ge 1$, and $0 \le r < k$.  Thus
\begin{equation}
	a_n \le (a_k)^p \, (a_1)^r
\end{equation}
by submultiplicativity.

	We can rewrite this as
\begin{equation}
	(a_n)^{1/n} \le (a_k)^{1/(k + (r/p))} \, (a_1)^{r/n}.
\end{equation}
Using this one can check that $\limsup_{n \to \infty} (a_n)^{1/n}$ is
less than or equal to $(a_k)^{1/k}$ for all positive integers $k$,
from which it follows that $\{(a_l)^{1/l}\}_{l=1}^\infty$ converges to
the infimum of the $a_j$'s, as desired.

	Let us review some notions which will be used here frequently.
If $V$ is a vector space over the real or complex numbers, then a
function $N(v)$ on $V$ is called a \emph{seminorm} if $N(v)$ is a
nonnegative real number for all $v \in V$ which is equal to $0$ when
$v = 0$ and satisfies
\begin{equation}
	N(\alpha \, v) = |\alpha| \, N(v)
\end{equation}
for all real or complex numbers $\alpha$, as appropriate, and all $v
\in V$, and
\begin{equation}
	N(v + w) \le N(v) + N(w)
\end{equation}
for all $v, w \in V$.  If furthermore $N(v) > 0$ when $v$ is a nonzero
vector in $V$, then $N$ is said to be a \emph{norm} on $V$.

	Let $\mathcal{A}$ be a real or complex vector space.  We say
that $\mathcal{A}$ is an \emph{algebra} if $\mathcal{A}$ is also
equipped with a binary operation of multiplication
\begin{equation}
	(x, y) \mapsto x \, y
\end{equation}
for $x, y \in \mathcal{A}$ which is associative and linear in each
variable $x$, $y$.

	Suppose that $\mathcal{A}$ is a real or complex algebra and
that $\|\cdot \|$ is a norm on $\mathcal{A}$.  We say that
$\mathcal{A}$ is a \emph{normed algebra} if we also have that
\begin{equation}
	\|x \, y\| \le \|x\| \, \|y\|
\end{equation}
for all $x, y \in \mathcal{A}$.

	Let $\mathcal{A}$ be a real or complex normed algebra with
norm $\|\cdot \|$, and let $x$ be an element of $\mathcal{A}$.  For
each positive integer $j$, let $x^j$ be the element of $\mathcal{A}$
which is a product of $j$ $x$'s.  The sequence of nonnegative real
numbers $a_j = \|x^j\|$ is submultiplicative, and therefore
$\{\|x^j\|^{1/j}\}_{j=1}^\infty$ converges as a sequence of real
numbers by the earlier discussion.

	Let $V$ be a real or complex vector space.  The linear
transformations from $V$ to itself form an algebra, using composition
of mappings as multiplication.

	Suppose further that $V$ is equipped with a norm $\|\cdot \|$.
A linear transformation $T$ on $V$ is said to be \emph{bounded} if
there is a nonnegative real number $k$ such that
\begin{equation}
	\|T(v)\| \le k \, \|v\|
\end{equation}
for all $v \in V$.  In this event we can define the operator norm
$\|T\|_{op}$ of $T$ by
\begin{equation}
	\|T\|_{op} = \sup \{\|T(v)\| : v \in V, \ \|v\| \le 1\},
\end{equation}
which is the same as the optimal value of $k$ for the preceding
inequality.

	The bounded linear operators on a normed vector space also
form an algebra.  The operator norm defines a norm in such a way that
the algebra of bounded linear operators becomes a normed algebra.

	Let $V$ be a real or complex vector space with norm $\| \cdot
\|$, and observe that
\begin{equation}
	d(v, w) = \|v - w\|
\end{equation}
defines a metric on $V$.  If $V$ is complete with respect to this
metric, in the sense that every Cauchy sequence in $V$ converges to
some element of $V$, then $V$ is called a \emph{Banach space}.

	In any normed vector space $(V, \|\cdot \|)$, a series
$\sum_{j=1}^\infty v_j$ with terms in $V$ is said to converge if the
sequence of partial sums converges.  The series $\sum_{j=1}^\infty
v_j$ converges absolutely if $\sum_{j=1}^\infty \|v_j\|$ converges as
a series of nonnegative real numbers.  Absolute convergence implies
that the sequence of partial sums forms a Cauchy sequence.

	Thus in a Banach space every absolutely convergent infinite
series converges.  Conversely one can show that if every absolutely
convergent series in a normed vector space converges, then the space
is complete.

	From now on in these notes we assume that a real or complex
algebra $\mathcal{A}$ is equipped with a nonzero multiplicative
identity element $e$.  For linear operators on a vector space $V$, the
identity transformation $I$ on $V$ plays this role, since the
composition of any linear transformation $T$ on $V$ with the identity
transformation $I$ is equal to $T$.

	We also assume that in a normed algebra $\mathcal{A}$ the norm
of the multiplicative identity element $e$ is equal to $1$.  If $V$ is
a normed vector space of positive dimension, then the algebra of
bounded linear operators on $V$ equipped with the operator norm has
this property, since the operator norm of the identity operator is
automatically equal to $1$.

	By a \emph{Banach algebra} we mean a normed algebra which is a
Banach space as a normed vector space.  As in the preceding paragraphs
we assume that the algebra contains a nonzero multiplicative identity
element whose norm is equal to $1$.  If $V$ is a Banach space, then
the algebra of bounded linear operators on $V$ is a Banach algebra
with respect to the operator norm.

	Let $\mathcal{A}$ be a real or complex Banach algebra with
multiplicative identity element $e$ and norm $\|\cdot \|$.  For each
nonnegative integer $n$ we have that
\begin{equation}
	(e - x) \bigg(\sum_{j=0}^n x^j \bigg)
		= \bigg(\sum_{j=0}^n x^j \bigg) (e - x)
			= e - x^{n+1},
\end{equation}
where $x^j$ is interpreted as being equal to $e$ when $j = 0$.

	Suppose that $\sum_{j=0}^\infty x^j$ converges in
$\mathcal{A}$.  As usual, this implies that
\begin{equation}
	\lim_{m \to \infty} x^m = 0
\end{equation}
in $\mathcal{A}$.  It follows that $e - x$ is invertible in
$\mathcal{A}$ in this case, and that its inverse is equal to
$\sum_{j=0}^\infty x^j$.

	If $\|x\| < 1$, or if $\|x^k\| < 1$ for any positive integer
$k$, then $\sum_{j=0}^\infty x^j$ converges absolutely, and hence
converges.  It was noted in class that convergence of
$\sum_{j=0}^\infty x^j$ implies that $\|x^k\| < 1$ for all
sufficiently large $k$, since $x^k \to 0$ as $k \to \infty$.

	Suppose that $x$, $y$ are elements of $\mathcal{A}$
with $x$ invertible and 
\begin{equation}
	\|y - x\| < \frac{1}{\|x^{-1}\|}.
\end{equation}
We can write $y$ as $x \, (e + x^{-1} (y - x))$ or as $(e + (y - x) \,
x^{-1}) \, x$ and conclude that $y$ is invertible.  Thus the
invertible elements of $\mathcal{A}$ form an open subset of
$\mathcal{A}$.

	If $x \in \mathcal{A}$, let $\rho(x)$ denote the resolvent set
associated to $x$, which is the set of real or complex numbers
$\lambda$, as appropriate, such that $\lambda \, e - x$ is invertible
in $\mathcal{A}$.  Let $\sigma(x)$ denote the spectrum of $x$, which
is the set of real or complex numbers $\lambda$, as appropriate, such
that $\lambda \, e - x$ is not invertible in $\mathcal{A}$.

	Let $x$ be an element of $\mathcal{A}$, and let $\lambda$ be a
real or complex number, as appropriate.  If $\lambda^n \, e - x^n$ is
invertible in $\mathcal{A}$, then so is $\lambda \, e - x$, which is
the same as saying that $\lambda^n \in \rho(x^n)$ implies that
$\lambda \in \rho(x)$.  Hence $\lambda \in \sigma(x)$ implies that
$\lambda^n \in \sigma(x^n)$.

	If $|\lambda| > \|x\|$, then $\lambda \, e - x$ is invertible,
and $\lambda \in \sigma(x)$.  For each positive integer $n$, if
$|\lambda|^n > \|x^n\|$, then $\lambda \in \rho(x)$.

	It follows from the earlier discussion that the resolvent set
is always an open set of real or complex numbers, as appropriate.  The
spectrum is therefore a closed set, which is in fact compact, because
it is bounded.

	In the complex case a famous result states that the spectrum
of an element $x$ of $\mathcal{A}$ is always nonempty.  For if
$\sigma(x)$ were empty, so that $\lambda \, e - x$ is invertible for
all complex numbers $\lambda$, then $(\lambda \, e - x)^{-1}$ would
define a nonzero holomorphic $\mathcal{A}$-valued function on the
complex plane which tends to $0$ as $\lambda \to 0$.  An extension of
Liouville's theorem would imply that this function is equal to $0$ for
all $\lambda$, a contradiction.

	In either the real or complex case, if $x \in \mathcal{A}$ and
$\sigma(x) \ne \emptyset$, let $r(x)$ denote the spectral radius of
$x$, which is the maximum of $|\lambda|$ over all $\lambda \in
\sigma(x)$.  Thus $r(x) \le \|x^n\|^{1/n}$ for all positive integers
$n$.

	Let $R(x)$ be equal to the infimum of $\|x^n\|^{1/n}$ over all
positive integers $n$, which is the same as the limit of
$\|x^n\|^{1/n}$ as $n \to \infty$, as we have seen.  In the complex
case a famous result states that $r(x) = R(x)$.  The idea is that $(e
- \alpha \, x)^{-1}$ defines a holomorphic $\mathcal{A}$-valued
function of $\alpha$ on the disk $|\alpha| < 1/r(x)$, and that the
series $\sum_{j=0}^\infty \alpha^j \, x^j$ consequently converges
everywhere on this disk.

	Suppose that $\phi$ is a nonzero homomorphism from
$\mathcal{A}$ into the real or complex numbers, as appropriate.  This
means that $\phi$ is a linear mapping from $\mathcal{A}$ into the real
or complex numbers, and that $\phi(x \, y) = \phi(x) \, \phi(y)$ for
all $x, y \in \mathcal{A}$.

	In particular, $\phi(e) = 1$.  More generally, if $x$ is an
invertible element of $\mathcal{A}$, then $\phi(x) \ne 0$.
If $\phi(x) = 0$ for some $x \in \mathcal{A}$, then $x$ is not
invertible.

	If $x \in \mathcal{A}$, then $\phi(x) \in \sigma(x)$.  This is
because $\phi$ applied to $\phi(x) \, e - x$ is equal to $0$.

	As a result,
\begin{equation}
	|\phi(x)| \le R(x) \le \|x\|
\end{equation}
for all $x \in \mathcal{A}$.  This implies that $\phi$ is continuous
as a mapping from $\mathcal{A}$ into the real or complex numbers.

	In the complex case for a commutative Banach algebra there are
abstract arguments to the effect that every complex number in the
spectrum of an element of the Banach algebra occurs as the value of a
nonzero homomorphism from the algebra into the complex numbers at the
element of the Banach algebra.  This implies that $R(x)$ is a seminorm
on $\mathcal{A}$.  Let us check this directly for commutative normed
algebras in general.

	Of course $R(x)$ is always a nonnegative real number, and it
is easy to see that $R(\alpha \, x) = |\alpha| \, R(x)$ for all real
or complex numbers $\alpha$, as appropriate, and all $x \in
\mathcal{A}$.  To show that $R(x)$ defines a seminorm on
$\mathcal{A}$, it remains to check the triangle inequality.

	Let $\{a_j\}_{j=0}^\infty$, $\{b_k\}_{k=0}^\infty$ be two
submultiplicative sequences of nonnegative real numbers, with $a_0 =
b_0 = 1$.  Define $c_n$ for each nonnegative integer $n$ by
\begin{equation}
	c_n = \sum_{j=0}^n {n \choose j} \, a_j \, b_{n-j},
\end{equation}
so that $c_0 = 1$.

	Actually, $\{c_n\}_{n=0}^\infty$ is a submultiplicative
sequence of nonnegative real numbers in this situation.  Explicitly,
\begin{equation}
	\sum_{j=0}^{m+n} {m + n \choose j} \, a_j \, b_{m + n - j}
		\le \bigg(\sum_{k=0}^m {m \choose k} \, a_k \, b_{m - k} \bigg)
	\cdot \bigg(\sum_{l = 0}^n {n \choose l} \, a_l \, b_{n - l} \bigg).
\end{equation}
This is not difficult to verify, multiplying out the product on the
right and using identities from the binomial theorem.

	Thus $\lim_{n \to \infty} (c_n)^{1/n}$ exists.  The next step
is that
\begin{equation}
	\lim_{n \to \infty} (c_n)^{1/n} 
  \le \lim_{j \to \infty} (a_j)^{1/j} + \lim_{k \to \infty} (b_k)^{1/k},
\end{equation}
which is not too difficult to show.

	In connection with this, notice that if $t_1, \ldots, t_m$
are nonnegative real numbers, then
\begin{equation}
	\max (t_1, \ldots, t_m) \le t_1 + \cdots + t_m
		\le m \cdot \max (t_1, \ldots, t_m).
\end{equation}
If $\{t_{p, n}\}_{n=0}^\infty$, $p = 1, \ldots, m$, are sequences
of nonnegative real numbers, then
\begin{equation}
	\limsup_{n \to \infty} \max (t_{1, n}, \ldots, t_{m, n})^{1/n}
		= \limsup_{n \to \infty} (t_{1, n} + \cdots + t_{m, n})^{1/n}.
\end{equation}

	Suppose that $x, y \in \mathcal{A}$, and put $a_j = \|x^j\|$,
$b_k = \|y^k\|$.  If $c_n$ is as before, then
\begin{equation}
	\|(x + y)^n\| \le c_n.
\end{equation}
It follows that
\begin{equation}
	R(x + y) \le R(x) + R(y),
\end{equation}
which is what we wanted.

	One can also check that $R(x)$ is compatible with
multiplication in the same way as for a normed algebra.  Namely $R(e)
= 1$ and $R(x \, y) \le R(x) \, R(y)$ when $x, y \in \mathcal{A}$ and
$\mathcal{A}$ is a commutative normed algebra.

	Let us look at some special cases.  Let $C({\bf T})$ denote
the algebra of continuous complex-valued functions on the unit circle
${\bf T}$ in the complex plane, equipped with the supremum norm,
$\|f\| = \sup \{|f(z)| : z \in {\bf T}\}$.  In this case $\|f^n\| =
\|f\|^n$ for every positive integer $n$.

	Now consider the continuous functions $f(z)$ on the unit
circle ${\bf T}$ which can be expressed as $\sum_{j=0}^\infty a_j \,
z^j + \sum_{l=1}^\infty a_{-l} \, \overline{z}^l$, where
$\overline{z}$ denotes the complex conjugate of the complex number $z$
and where the coefficients $a_j$ are complex numbers such that
$\sum_{j=-\infty}^\infty |a_j|$ converges.  We define a norm $\|f\|_1$
to be equal to this sum.

	This expansion for $f$ is simply the Fourier series of $f$,
and the point is that we restrict our attention to continuous
functions whose Fourier series is absolutely summable.  We certainly
have that $\|f\| \le \|f\|_1$ for such a function, which is to say
that the supremum of $|f(z)|$ is bounded by $\sum_{j=-\infty}^\infty
|a_j|$.  One can check that if $f_1$, $f_2$ are two functions on the
unit circle with absolutely summable Fourier series, then so is the
product $f_1 \, f_2$ and $\|f_1 \, f_2\|_1 \le \|f_1\|_1 \,
\|f_2\|_1$.

	If $f$ is a function on the unit circle with absolutely
summable Fourier series, then a famous result states that
\begin{equation}
	\lim_{n \to \infty} (\|f^n\|_1)^{1/n} = \|f\|.
\end{equation}
Because $R(f) = \lim_{n \to \infty} (\|f^n\|_1)^{1/n}$ defines a
seminorm, it is sufficient to check that $R(f) = \|f\|$ for a dense
class of functions such as trigonometric polynomials.

	Now let $V$ be the vector space of real or complex sequences
$x = \{x_j\}_{j=1}^\infty$, as one might prefer, in which at most
finitely many terms are nonzero.  If $p$ is a real number, $1 \le p <
\infty$, then we can define a norm $\|x\|_p$ on $V$ by
\begin{equation}
	\|x\|_p = \bigg(\sum_{j=1}^\infty |x_j|^p \bigg)^{1/p}.
\end{equation}
For $p = \infty$ put
\begin{equation}
	\|x\|_\infty = \max \{|x_j| : j \ge 1\}.
\end{equation}

	Let $\{\alpha_j\}_{j=1}^\infty$ be a monotone decreasing
sequence of nonnegative real numbers.  Define a linear transformation
$T$ on $V$ by saying that the $j$th term of the sequence $T(x)$ is
equal to $\alpha_j \, x_{j+1}$.

	Define a sequence $\{a_l\}_{l=1}^\infty$ of nonnegative real
numbers by
\begin{equation}
	a_l = \prod_{j=1}^j \alpha_j.
\end{equation}
For each positive integer $l$ and each $p$, $1 \le p \le \infty$,
$a_l$ is the operator norm of $T^l$ with respect to the norm $\|x\|_p$
on $V$.  One can check directly that $\{a_l\}_{l=1}^\infty$ is a
submultiplicative sequence in this case, and that $\lim_{l \to \infty}
(a_l)^{1/l}$ is equal to $\lim_{j \to \infty} \alpha_j$, which is the
same as the infimum of the $\alpha_j$'s.

\end{document}